\documentclass{amsart}

\usepackage[all]{xy}
\usepackage{amscd}
\usepackage[mathscr]{eucal}
\usepackage{amssymb}

\newcommand{\colim}{\operatornamewithlimits{colim}}
\newcommand{\hocolim}{\operatornamewithlimits{hocolim}}
\newcommand{\holim}{\operatornamewithlimits{holim}}
\DeclareMathOperator{\Map}{Map}
\newtheorem{lem}{Lemma}[section]

\newtheorem{prop}[lem]{Proposition}

\newtheorem{theo}[lem]{Theorem}
\newtheorem{coro}[lem]{Corollary}

\newtheorem*{theorem}{Theorem}

\begin{document}
\title{Split Injectivity of the Baum-Connes Assembly Map}
\author{David Rosenthal}
\address{Department of Mathematics and Statistics, McMaster University, Hamilton, ON L8S 4K1, Canada}
\email{rosend@math.mcmaster.ca}

\begin{abstract}
	In this work, the continuously controlled techniques developed by Carlsson and Pedersen are used to prove that the Baum-Connes map is a split injection for groups satisfying certain geometric conditions.
\end{abstract}

\maketitle

\section{Introduction}

	In~\cite{higson}, Higson used analytic techniques to prove that the Baum-Connes map
\[ KK_i^{\Gamma} (C_0(E\Gamma(\mathfrak{f}));\mathbb{C}) \to K_i(C^*_r\Gamma), \]
where $\Gamma$ is a discrete group and $E\Gamma(\mathfrak{f})$ denotes the universal space for $\Gamma$-actions with finite isotropy, is injective when $E\Gamma(\mathfrak{f})$ satisfies certain geometric conditions. The goal of this paper is to use continuously controlled algebra to prove a strengthening of this result, namely that under these conditions, the Baum-Connes map is a {\it split} injection. This tells us that the left-hand side of the map is a direct summand of the right-hand side. Thus, a piece of $K_i(C^*_r\Gamma)$ is understood since there are tools available for calculating homology groups. The geometric assumptions placed on $E\Gamma(\mathfrak{f})$ are satisfied by a large class of groups, including virtually polycyclic groups and word hyperbolic groups.

	The following theorem is proven in this paper.

\begin{theorem}
	Let $\Gamma$ be a discrete group, and let $\mathcal{E}=E\Gamma(\mathfrak{f})$ be the universal space for $\Gamma$-actions with finite isotropy, where $\mathfrak{f}$ denotes the family of finite subgroups of $\Gamma$. Assume that $\mathcal{E}$ is a finite $\Gamma$-CW complex admitting a compactification, $X$, (i.e., $X$ is compact and $\mathcal{E}$ is an open dense subset) such that
	\begin{enumerate}
		\item[1.] the $\Gamma$-action extends to $X$;
		\item[2.] $X$ is metrizable;
		\item[3.] $X^G$ is contractible for every $G\in \mathfrak{f}$;
		\item[4.] ${\mathcal{E}}^G$ is dense in $X^G$ for every $G\in \mathfrak{f}$;
		\item[5.] compact subsets of $\mathcal{E}$ become small near $Y=X-\mathcal{E}$.  That is, for every compact subset $K \subset \mathcal{E}$ and for every neighborhood $U \subset X$ of $y \in Y$, there exists a neighborhood $V \subset X$ of $y$ such that $g \in \Gamma$ and $gK \cap V \neq \emptyset$ implies $gK \subset U$.
	\end{enumerate}
Then the Baum-Connes map, $KK_i^{\Gamma}(C_0(\mathcal{E});\mathbb{C}) \to K_i(C^*_r\Gamma)$, is a split injection.
\end{theorem}

	The proof uses continuously controlled topology, based on Carlsson and Pedersen's work in~\cite{cp} where they split the assembly map in algebraic $K$- and $L$-theory for torsion free groups satisfying the same geometric conditions as are used here. Their results were generalized to groups with torsion in~\cite{me2,me3}. By interpreting the Baum-Connes map as an assembly map, continuously controlled algebra can enter the picture. In the algebraic $K$- and $L$-theory cases, the assembly map was split by realizing it as a map of fixed spectra. This allowed homotopy fixed sets to be used. More precisely, if $S\to T$ is a $\Gamma$-equivariant map of spectra, then the $\Gamma$-equivariant map $E\Gamma(\mathfrak{f})\to \bullet$, induces the diagram
\[ \xymatrix{
	S^{\Gamma} \ar[d] \ar[r] & T^{\Gamma} \ar[d] \\
	S^{h\mathfrak{f} \Gamma} \ar[r]^{h} & T^{h\mathfrak{f} \Gamma} } \]
where $S^{h\mathfrak{f} \Gamma}$ denotes the set of all $\Gamma$-equivariant maps from $E\Gamma(\mathfrak{f})$ to $S$. The homotopy fixed set, $S^{h\mathfrak{f} \Gamma}$, behaves well with respect to homotopy. In~\cite{cp, me2, me3}, the splittings were obtained by showing that the leftmost vertical map and the bottom map in the diagram were weak homotopy equivalences.

	In order to enter the world of topological $K$-theory, additional structure must be added to the continuously controlled categories. In the process, the hope for realizing the Baum-Connes map as a map of fixed sets is lost. However, it is closely related to a map of fixed sets. Homotopy fixed sets can still be used to split the Baum-Connes map. The added structure also means that work must be done in order to extend the continuously controlled results from~\cite{cp, me2}.
	
	Section 2 provides the background for the paper, including the basic definitions from continuously controlled algebra, a brief discussion of the $\mathbb{K}^{Top}$-functor, the continuously controlled formulation of the Baum-Connes assembly map, and the necessary facts about homotopy fixed sets. In Section 3, the groundwork for the proof of the main theorem is laid. In this section, the continuously controlled algebra results from~\cite{cp} that we require are extended to this setting. It turns out that in the proof of the main theorem, it will be necessary to know that the reduced Steenrod homology (described in this section) of the quotient of a contractible compact metrizable space by a finite group is trivial. This is proven making use of the Conner Conjecture. In Section 4, the main theorem is proven. The most difficult step in the proof is handled by filtering the space $X$ by {\it conjugacy classes of fixed sets}.
	
	I would like to thank my uncle Peter Rosenthal, for many helpful discussions on operator theory.

\section{Preliminaries}

\subsection{Continuously Controlled Algebra}

	Let $\Gamma$ be a discrete group, $X$ a $\Gamma$-space, and $Y$ a closed $\Gamma$-invariant subspace of $X$. Let $E=X-Y$. The  complex vector space with basis $E \times \Gamma \times \mathbb{N}$ is denoted $\mathbb{C}[E \times \Gamma]^{\infty}$. If $H$ is a sub-vector space of $\mathbb{C}[E \times \Gamma]^{\infty}$, denote $H \cap \mathbb{C}[x \times \Gamma]^{\infty}$ by $H_x$, where $x \in E$. The {\it continuously controlled category} $\mathcal{B}(X,Y;\mathbb{C})$, has objects $H$, where
\begin{enumerate}
	\item[(i)] $H=\bigoplus_{x \in E} H_x$;
	\item[(ii)] $H_x$ is a finite dimensional complex vector space with basis contained in $\{ x \} \times \Gamma \times \mathbb{N}$;
	\item[(iii)] $\{x \in E \, | \, H_x \neq 0 \}$ is locally finite in $E$.
\end{enumerate}
Morphisms are all linear operators $\phi:H\to K$, that are {\it continuously controlled}. This means that for every $y\in Y$ and every neighborhood $U \subseteq X$ of $y$, there exists a neighborhood $V \subseteq X$ of $y$ such that the components of $\phi$, $\phi^x_z:H_x \to K_z$ and $\phi^z_x: K_z \to H_x$, are zero whenever $x \in V$ and $z \notin U$.

	The definition of this category can also be made for an arbitrary ring. Although such generality is not needed in this paper, we will encounter situations in which the complex group ring $\mathbb{C} G$ is used, where $G$ is a subgroup of $\Gamma$. The only changes in the definition are that the objects of $\mathcal{B}(X,Y;\mathbb{C} G)$ are free $\mathbb{C} G$-modules, where each $H_x$ is a finite dimensional free $\mathbb{C} G$-module, and morphisms have the additional feature that their components are $\mathbb{C}G$-homomorphisms.
	
	Note that $\mathbb{C}[E \times \Gamma]^{\infty}$ comes equipped with a free $\Gamma$-action. This action puts an interesting action of $\Gamma$ on $\mathcal{B}(X,Y;\mathbb{C})$. If we assume that $\Gamma$ acts on $E$ with finite isotropy, then the fixed category, $\mathcal{B}^{\Gamma}(X,Y;\mathbb{C})$, has those objects, $H$, in $\mathcal{B}(X,Y;\mathbb{C})$ that satisfy the conditions
\begin{enumerate}
  \item[1.] $H_{\gamma x}\cong H_x$ for every $\gamma \in \Gamma$, and
  \item[2.] $H_x$ is a finitely generated free $\mathbb{C}{\Gamma}_x$-module. 
\end{enumerate}  
This implies that $\bigoplus_{x'\in [x]} H_{x'}$ is a finitely generated free $\mathbb{C}\Gamma$-module, where $[x]=\{ \gamma x \, | \, \gamma \in \Gamma \}$. The morphisms in $\mathcal{B}^{\Gamma}(X,Y;\mathbb{C})$ are those morphisms, $\phi$, in $\mathcal{B}(X,Y;\mathbb{C})$ that are $\Gamma$-equivariant, i.e., $\gamma {\phi}^x_y {\gamma}^{-1}=\phi^{\gamma x}_{\gamma y}$ for all $\gamma \in \Gamma$ and all $x,y \in E$. Note that if a finite group $H$ acts trivially on $E$ then the fixed category $\mathcal{B}^H(X,Y;\mathbb{C})$ is equivalent to $\mathcal{B}(X,Y;\mathbb{C} H)$. Also note that the definition of the continuously controlled category used here is not designed to handle infinite isotropy. If there are infinite isotropy subgroups, then every object in the fixed category will be the zero object. However, for applications to the Baum-Connes conjecture only the case of finite isotropy is needed. 
	
	In order to discuss topological $K$-theory, our categories need to be modified. The category $C^*_r\mathcal{B}(X,Y;\mathbb{C})$, is defined to have the same objects as $\mathcal{B}(X,Y;\mathbb{C})$, but morphism sets are defined as follows. Every object, $H$, has a unique Hilbert space completion, $\mathcal{H}$, and every bounded morphism, $T:H \to K$, extends uniquely to a bounded operator $T:\mathcal{H} \to \mathcal{K}$. The set of morphisms from $H$ to $K$ in $C^*_r\mathcal{B}(X,Y;\mathbb{C})$ is obtained by first considering the subgroup of continuously controlled morphisms from $H$ to $K$ that are bounded operators, and then taking the closure of their extensions inside  $\mathfrak{B}(\mathcal{H},\mathcal{K})$. Notice that by doing this, strict control over morphisms in $C^*_r\mathcal{B}(X,Y;\mathbb{C})$ is lost. However, any morphism in $C^*_r\mathcal{B}(X,Y;\mathbb{C})$ can be approximated arbitrarily closely by a controlled morphism in $\mathcal{B}(X,Y;\mathbb{C})$. Since group actions are continuous, we similarly define $C^*_r\mathcal{B}^{\Gamma}(X,Y;\mathbb{C})$. At first glance this category appears to be the fixed category $(C^*_r\mathcal{B}(X,Y;\mathbb{C}))^{\Gamma}$, but this is not the case in general. A morphism in $C^*_r\mathcal{B}^{\Gamma}(X,Y;\mathbb{C})$ is a limit of $\Gamma$-equivariant continuously controlled morphisms, whereas in $(C^*_r\mathcal{B}(X,Y;\mathbb{C}))^{\Gamma}$ there can be $\Gamma$-equivariant morphisms that are limits of non-equivariant continuously controlled morphisms. There is, however, an inclusion functor from $C^*_r\mathcal{B}^{\Gamma}(X,Y;\mathbb{C})$ to $(C^*_r\mathcal{B}(X,Y;\mathbb{C}))^{\Gamma}$. If $\Gamma$ is a finite group, then this functor is an equivalence.
	
\begin{lem}\label{finite}
	If $\Gamma$ is a finite group, then $C^*_r\mathcal{B}^{\Gamma}(X,Y;\mathbb{C}) \cong (C^*_r\mathcal{B}(X,Y;\mathbb{C}))^{\Gamma}$.
\end{lem}

\begin{proof}
	We need only resolve the issue mentioned above. Let $T$ be a morphism in $(C^*_r\mathcal{B}(X,Y;\mathbb{C}))^{\Gamma}$. Then $T$ is a limit of continuously controlled bounded operators $\{ T_n \}$. Let
\[ S_n=\dfrac{1}{|\Gamma |}\sum_{g \in \Gamma} gT_ng^{-1}. \]
Then $S_n$ is a $\Gamma$-equivariant continuously controlled bounded operator, and since $T$ is $\Gamma$-equivariant, $\{ S_n \}$ converges to $T$.
\end{proof}
	
	Two controlled categories that are of particular interest to us are $\mathcal{B}(CX,CY \cup X,p_X;\mathbb{C})$ and $\mathcal{B}(\Sigma X,\Sigma Y,p_X;\mathbb{C})$, where $CX$ denotes the cone of $X$, $\Sigma X$ denotes the unreduced suspension of $X$, and $p_X: X\times (0,1) \to X$ is the projection map. These categories have the same objects as $\mathcal{B}(CX,CY \cup X;\mathbb{C})$ and $\mathcal{B}(\Sigma X,\Sigma Y;\mathbb{C})$ respectively, but their control conditions on morphisms differ along $Y \times (0,1)$, where they are only required to be $p_X$-{\it controlled}. This means that for every $(y,t)\in Y \times (0,1)$ and every neighborhood $U \subseteq X$ of $y$, there is a neighborhood $V$ of $(y,t)$ such that $\phi^a_b =0$ and $\phi^b_a =0$ whenever $a \in V \cap p_X^{-1}(U)$ and $b \notin p_X^{-1}(U)$.

	In~\cite{karoubi}, Karoubi introduced the notion of an $\mathcal{A}$-{\it filtered} additive category, where $\mathcal{A}$ is a full subcategory of $\mathcal{U}$. The associated {\it quotient category}, $\mathcal{U} / \mathcal{A}$, has the same objects as $\mathcal{U}$, but two morphisms are identified if their difference factors through $\mathcal{A}$. The {\it support at infinity} of an object $H$ in $\mathcal{B}(X,Y;\mathbb{C})$, denoted ${\rm supp}_{\infty}(H)$, is the set of limit points of $\{x \, | \, H_x \neq 0\}$. If $C$ is a closed  subspace of $Y$, then the category $\mathcal{B}(X,Y;\mathbb{C})_C$ is the full subcategory of $\mathcal{B}(X,Y;\mathbb{C})$ on objects $H$, with ${\rm supp}_{\infty}(H)\subseteq C$. The category $\mathcal{B}(X,Y;\mathbb{C})$ is $\mathcal{B}(X,Y;\mathbb{C})_C$-filtered. If $W$ is an open subspace of $Y$, then the {\it germ category} $\mathcal{B}(X,Y;\mathbb{C})^W$ has the same objects as $\mathcal{B}(X,Y;\mathbb{C})$, but morphisms are identified if they agree in a neighborhood of $W$. It is equivalent to the quotient category $\mathcal{B}(X,Y;\mathbb{C}) / \mathcal{B}(X,Y;\mathbb{C})_C$, when $C=Y-W$~\cite{cp}. The category $C^*_r\mathcal{B}(X,Y;\mathbb{C})_C$ is the full subcategory of $C^*_r\mathcal{B}(X,Y;\mathbb{C})$ on objects $H$, with ${\rm supp}_{\infty}(H)\subseteq C$, and $C^*_r\mathcal{B}(X,Y;\mathbb{C})$ is $C^*_r\mathcal{B}(X,Y;\mathbb{C})_C$-filtered. Two morphisms in the corresponding quotient category are identified if their difference can be approximated arbitrarily closely by continuously controlled bounded operators that factor through objects in $\mathcal{B}(X,Y;\mathbb{C})_C$. This is equivalent to saying that two morphisms are identified if their difference can be approximated arbitrarily closely by continuously controlled bounded operators that are zero in a neighborhood of $W=Y-C$. Thus, we will denote this quotient category by $C^*_r\mathcal{B}(X,Y;\mathbb{C})^{Y-C}$. The reason for studying these {\it Karoubi filtrations} is that the sequence
\[ C^*_r\mathcal{B}(X,Y;\mathbb{C})_C \to C^*_r\mathcal{B}(X,Y;\mathbb{C}) \to C^*_r\mathcal{B}(X,Y;\mathbb{C})^{Y-C} \]
induces a fibration of spectra after applying $\mathbb{K}^{Top}$. Such sequences therefore become the main tool in this theory. The above discussion also works for fixed categories if we use $\Gamma$-invariant subspaces.

\subsection{The $\mathbb{K}^{Top}$ Functor}

	There are two functors called $\mathbb{K}^{Top}$, one from the category of $C^*$-algebras to the category of spectra and the other from the category of $C^*$-categories to the category of spectra. The two functors are closely related. The definition of $\mathbb{K}^{Top}$ of a $C^*$-algebra is simple to state. For a given $C^*$-algebra $A$, $\mathbb{K}^{Top}(A)=\Omega^{-1}\mathbb{U}(A)$, where $\mathbb{U}(A)$ is the $\Omega$-spectrum whose first term is the infinite loop space $U(A)=\{ u\in A\otimes \mathfrak{K} \, | \, 1+u {\rm \, is \, a \, unitary \, in \, the \, unitization \, of \,} A\otimes \mathfrak{K} \}$, which satisfies Bott periodicity~\cite{higroeped}. The definition of $\mathbb{K}^{Top}$ of a $C^*$-category is not as easy to state and will not be presented here (an excellent reference is~\cite{hp}). However, the catch phrase is that the $K$-theory of a $C^*$-category is the direct limit of the $K$-theory of the endomorphisms of the objects of the category. For example, $\mathbb{K}^{Top}(C^*_r\mathcal{B}(\bullet, \emptyset; \mathbb{C} \Gamma))\simeq \mathbb{K}^{Top}(C^*_r(\Gamma))$.
	
	In~\cite{higroeped}, the $K$-theory of continuously controlled $C^*$-categories was identified with a certain $C^*$-algebra constructed by John Roe. This correspondence is used to establish certain properties of the $\mathbb{K}^{Top}$-functor that are necessary for using the continuously controlled techniques. Under this identification, a filtered subcategory corresponds to a closed two-sided ideal in the Roe algebra. Since a short exact sequence of $C^*$-algebras induces a fibration of spectra after applying $\mathbb{K}^{Top}$~\cite[4.4]{higroeped}, the Karoubi filtration
\[ C^*_r\mathcal{B}(X,Y;\mathbb{C})_C \to C^*_r\mathcal{B}(X,Y;\mathbb{C}) \to C^*_r\mathcal{B}(X,Y;\mathbb{C})^{Y-C} \]
also yields a fibration of spectra (up to homotopy) after applying $\mathbb{K}^{Top}$.

	It follows from the definition that $\mathbb{K}^{Top}$ of a countable $C^*$-direct sum of $C^*$-algebras $\bigoplus_i A_i$, is equivalent to the countable direct product $\prod_i \mathbb{K}^{Top}(A_i)$. This implies that $\mathbb{K}^{Top}$ of a countable $C^*$-direct sum of $C^*$-categories $\bigoplus_i \mathcal{A}_i$, is equivalent to the countable direct product $\prod_i \mathbb{K}^{Top}(\mathcal{A}_i)$. (A morphism $T$, in the $C^*$-category $\bigoplus_i \mathcal{A}_i$, is a sequence of uniformly bounded morphisms $\{ T_i \}$, where each $T_i$ is a morphism in $\mathcal{A}_i$.) The final property that we require is that $\mathbb{K}^{Top}$ of a fixed $C^*$-category $\mathcal{A}^{\Gamma}$, is equivalent to $\mathbb{K}^{Top}(\mathcal{A})^{\Gamma}$.

\subsection{The Baum-Connes Assembly Map}

	The Baum-Connes assembly map is the map of spectra
\[ \Omega\mathbb{K}^{Top}\Big(C^*_r\mathcal{B}^{\Gamma}(\mathcal{E} \times (0,1],\mathcal{E} \times 1;\mathbb{C})^{\mathcal{E} \times 1}\Big) \to \Omega\mathbb{K}^{Top}\Big(C^*_r\mathcal{B}(\bullet \times (0,1],\bullet \times 1;\mathbb{C}\Gamma)^{\bullet \times 1}\Big) \]
induced by collapsing $\mathcal{E}$ to a point, where $\mathcal{E}$ denotes the universal space for $\Gamma$-actions with finite isotropy. On homotopy groups this is the Baum-Connes map $KK_i^{\Gamma} (C_0(\mathcal{E});\mathbb{C}) \to K_i(C^*_r\Gamma)$~\cite[Theorem 7.6]{hp}. We would like to realize the Baum-Connes assembly map as a map of fixed spectra so that we can proceed to split it as in~\cite{cp, me2} (see the next section). However, this cannot be done. When $\mathcal{E}$ admits a compactification $X$, satisfying the conditions of the main theorem, the Baum-Connes assembly map is weakly homotopy equivalent to
\[ \Omega\mathbb{K}^{Top}(C^*_r\mathcal{B}^{\Gamma}(CX,CY \cup X,p_X;\mathbb{C})) \to \Omega\mathbb{K}^{Top}(C^*_r\mathcal{B}^{\Gamma}(\Sigma X,\Sigma Y,p_X;\mathbb{C})) \]
(see Section 3.1). Although this is not a map of fixed spectra, it comes close enough that it will still be possible to use the techniques developed for splitting those assembly maps that are maps of fixed spectra.

\subsection{Generalized Fixed Sets}

	In this section, we recall some facts about generalized homotopy fixed sets in the category of spectra. Let $S$ be a spectrum with $\Gamma$-action. The fixed set $S^{\Gamma}$, can be identified with the set $\Map_{\Gamma}(\bullet,S)$ of $\Gamma$-equivariant maps from a point into $S$. Let $\mathcal{F}$ be a family of subgroups of $\Gamma$. The homotopy fixed set associated to this family $S^{h \mathcal{F} \Gamma}$, is defined to be $\Map_{\Gamma}(E\Gamma(\mathcal{F}),S)$. Notice that the $\Gamma$-equivariant map $E\Gamma(\mathcal{F}) \to \bullet$
induces the following commutative diagram:
\[ \xymatrix{
	S^{\Gamma} \ar[d]_a \ar[r] & T^{\Gamma} \ar[d] \\
	S^{h\mathcal{F} \Gamma} \ar[r]^b & T^{h\mathcal{F} \Gamma}. } \]
If $a$ and $b$ are shown to be weak homotopy equivalences, then at the level of homotopy, the map $S^{\Gamma} \to T^{\Gamma}$ will have successfully been split. In~\cite{cp, me2, me3}, such a diagram was used to split the assembly maps in algebraic $K$- and $L$-theory by realizing them as maps of fixed sets. Although the Baum-Connes assembly map cannot be realized as a map of fixed sets, we will still be able to make use of the above diagram (see Section 4). The next two lemmas play a key role in such an approach. The author would like to thank Wolfgang L\"{u}ck for showing him the elegant proof of the following lemma.

\begin{lem}\label{lemmaa}
	Let $F:S \to T$ be an equivariant map of spectra with $\Gamma$-action. If $F^G:S^G \to T^G$ is a weak homotopy equivalence for every $G\in \mathcal{F}$, then $S^{h \mathcal{F} \Gamma} \simeq T^{h \mathcal{F} \Gamma}$.
\end{lem}

\begin{proof}
	Suppose for the moment that $\Map_{\Gamma}(X,S) \to \Map_{\Gamma}(X,T)$ is a weak homotopy equivalence for any finite dimensional $\Gamma$-CW complex $X$. Since the inclusion of the $n$-skeleton of $E\Gamma(\mathcal{F})$ into its $(n+1)$-skeleton $E\Gamma(\mathcal{F})_n \hookrightarrow E\Gamma(\mathcal{F})_{n+1}$, is a cofibration, $E\Gamma(\mathcal{F})=\colim_n E\Gamma(\mathcal{F})_n \simeq \hocolim_n E\Gamma(\mathcal{F})_n$. Therefore
\begin{align}
	&S^{h \mathcal{F} \Gamma} & &\simeq & & \Map_{\Gamma}\big(\hocolim_n E\Gamma(\mathcal{F})_n,S\big) & &\simeq & &\holim_n\Map_{\Gamma}(E\Gamma(\mathcal{F})_n,S) \notag\\
	& & & & & & &\simeq & &\holim_n\Map_{\Gamma}(E\Gamma(\mathcal{F})_n,T)\notag\\
	& & & & & & &\simeq & &T^{h \mathcal{F} \Gamma}.\notag
\end{align}	

	Let $X$ be a finite dimensional $\Gamma$-CW complex. We will prove that $\Map_{\Gamma}(X,S) \simeq \Map_{\Gamma}(X,T)$ by induction on the skeleta of $X$. 
\[ \Map_{\Gamma}(X_0,S)=\Map_{\Gamma}\Big(\bigsqcup_i \Gamma / H_i,S\Big)=\prod_i\Map_{\Gamma}(\Gamma / H_i,S)=\prod_i S^{H_i}, \]
where each $H_i \in \mathcal{F}$. Since $S^{H_i}\simeq T^{H_i}$ for every $i$, $\Map_{\Gamma}(X_0,S) \simeq \Map_{\Gamma}(X_0,T)$.

	Now assume that $\Map_{\Gamma}(X_{k-1},S) \simeq \Map_{\Gamma}(X_{k-1},T)$. Consider the following pushout diagram
\[ \xymatrix{
	\Map_{\Gamma}\big(\bigsqcup_j \Gamma / G_j \times S^{k-1},S\big) \ar[d] \ar[r] & \Map_{\Gamma}(X_{k-1},S) \ar[d] \\
	\Map_{\Gamma}\big(\bigsqcup_j \Gamma / G_j \times \mathbb{D}^k,S\big) \ar[r] & \Map_{\Gamma}(X_k,S) } \]
where each $G_j\in \mathcal{F}$. The equivariant map from $S$ to $T$ induces a commutative diagram relating the above pushout diagram to the corresponding one for $T$. Therefore the proof will be complete if we show that the maps $\Map_{\Gamma}(X_{k-1},S) \to \Map_{\Gamma}(X_{k-1},T)$ and $\Map_{\Gamma}\big(\bigsqcup_j \Gamma / G_j \times \mathbb{D}^k,S\big) \to \Map_{\Gamma}\big(\bigsqcup_j \Gamma / G_j \times \mathbb{D}^k,T\big)$ are weak homotopy equivalences. The first map is a weak homotopy equivalence by the induction hypothesis, and since
\[ \Map_{\Gamma}\big(\bigsqcup_j \Gamma / G_j \times \mathbb{D}^k,S\big)=\prod_j\Map_{\Gamma}(\Gamma / G_j \times \mathbb{D}^k,S) \simeq \prod_j\Map_{\Gamma}(\Gamma / G_j,S)= \prod_j S^{G_j} \]
the second map is also a weak homotopy equivalence.
\end{proof}

\begin{lem}\label{lemmab}
	Let $B$ be a $G$-spectrum, where $G\in \mathcal{F}$, and let $\Gamma$ act on $S=\prod_{\Gamma / G}B$ by identifying $S$ with $\Map_G(\Gamma,B)$ (in which $(\gamma f)(x)=f({\gamma}^{-1}x)$, where $\gamma \in \Gamma$). Then $S^{\Gamma} \simeq S^{h \mathcal{F} \Gamma}$.
\end{lem}

\begin{proof}
\[ S^{h \mathcal{F} \Gamma}= \Map_{\Gamma}(E\Gamma(\mathcal{F}),\Map_G(\Gamma,B))=\Map_G(E\Gamma(\mathcal{F}),B) \]
by evaluating at 1.
Since $E\Gamma(\mathcal{F})$ is $G$-equivariantly homotopy equivalent to $E\Gamma(\mathcal{F})^G$,
\begin{align}
	&\Map_G(E\Gamma(\mathcal{F}),B) & &\cong & &\Map_G(E\Gamma(\mathcal{F})^G,B) & &= & &\Map(E\Gamma(\mathcal{F})^G,B^G). \notag
\end{align}
Finally, $\Map(E\Gamma(\mathcal{F})^G,B^G)=B^G=S^{\Gamma}$ since $E\Gamma(\mathcal{F})^G$ is contractible.
\end{proof}

\section{Key Facts}

\subsection{$C^*$-algebras and Controlled Categories}

	In order to handle topological $K$-theory, additional structure needed to be placed on the continuously controlled categories. Recall that morphisms in $C^*_r\mathcal{B}^{\Gamma}(X,Y;\mathbb{C})$ are required to be bounded operators. Furthermore, morphism sets are completed, so strict control is lost. Therefore we must verify that the various results that we need from~\cite{cp, me2} can be extended to their $C^*_r$-counterparts. In general this is not so easy. For example, a functor between two controlled categories will extend to a functor in the $C^*_r$-setting if it sends bounded operators to bounded operators and Cauchy sequences of bounded operators to Cauchy sequences. However, checking whether or not a given functor possess these properties is non-trivial. Luckily, the various functors used in~\cite{cp,me2} are defined in such a way that they are guaranteed to extend to their $C^*_r$-analogues.
	
	A function $f:(X,Y) \to (X',Y')$ is {\it eventually continuous} if for every compact $K \subset X'-Y'$, $f^{-1}(K)$ has compact closure in $X-Y$, $f(X-Y) \subset X'-Y'$, and $f$ is continuous on $Y$. In \cite[Proposition 1.16]{cp}, it is proven that an eventually continuous function $f:(X,Y) \to (X',Y')$, induces a functor $U_f:\mathcal{B}(X,Y;\mathbb{C}) \to \mathcal{B}(X',Y';\mathbb{C})$ (well-defined up to natural equivalence), defined as follows. An object $H$ is sent to $U_f(H)$, where $(U_f(H))_{x'}=\bigoplus^{}_{x\in f^{-1}(x')} H_x$, and comes equipped with a unitary operator $U_H:H \to U_f(H)$, defined by $U_H(v)_{x'}=\bigoplus^{}_{x\in f^{-1}(x')} v_x$, where $v\in H$. If $T:H\to K$ is a continuously controlled morphism, then $U_f(T)=U_K T U_H^{-1}$. The unitary $U_H$ extends uniquely to a unitary, also denoted by $U_H$, on the Hilbert space completions of $H$ and $U_f(H)$. Now it is easy to see that $U_f$ extends to a functor $\mathscr{U}_f:C^*_r\mathcal{B}(X,Y;\mathbb{C}) \to C^*_r\mathcal{B}(X',Y';\mathbb{C})$, where $\mathscr{U}_f(H)=U_f(H)$, and $\mathscr{U}_f(T)=U_K T U_H^{-1}$ when $T$ is a morphism from $H$ to $K$ in $C^*_r\mathcal{B}(X,Y;\mathbb{C})$.
	
	The remainder of this section is devoted to proving the $C^*_r$-versions of the results in~\cite{cp} that we will need, making frequent use of {\it unitary functors} like the ones induced by eventually continuous maps. Throughout this section assume that $X$ is a compact metrizable space, $Y$ is a closed nowhere dense subset, and let $E=X-Y$.

\begin{lem}\label{equal}
	If $f_1$, $f_2:(X,Y) \to (X',Y')$ are eventually continuous maps and $f_1|_Y=f_2|_Y$, then $\mathscr{U}_{f_1}$ and $\mathscr{U}_{f_2}$ are naturally equivalent functors.
\end{lem}

\begin{proof}
	Let $H$ be an object in $C^*_r\mathcal{B}(X,Y;\mathbb{C})$. Let $U_H:H \to \mathscr{U}_{f_1}(H)$ be the unitary operator associated to $H$ by $\mathscr{U}_{f_1}$, and let $V_H:H \to \mathscr{U}_{f_2}(H)$ be the unitary operator associated to $H$ by $\mathscr{U}_{f_2}$. Then the natural equivalence, $\eta:\mathscr{U}_{f_1}\to \mathscr{U}_{f_2}$, is given by $\eta(H)=V_H \circ U_H^{-1}$.
\end{proof}

\begin{prop}\label{cone}
	There is a natural equivalence
\[ \mathscr{U}_f:C^*_r\mathcal{B}(X,Y;\mathbb{C}) \to C^*_r\mathcal{B}(CY,Y;\mathbb{C}). \]
That is, if $X'$ is a compact metrizable space, $Y'$ is a closed nowhere dense subset, and $h:(X,Y) \to (X',Y')$ is an eventually continuous map, then 
	\[ \xymatrix{
	C^*_r\mathcal{B}(X,Y;\mathbb{C}) \ar[d] \ar[r] & C^*_r\mathcal{B}(CY,Y;\mathbb{C}) \ar[d] \\
	C^*_r\mathcal{B}(X',Y';\mathbb{C}) \ar[r] & C^*_r\mathcal{B}(CY',Y';\mathbb{C}) } \]
commutes up to natural equivalence.
\end{prop}

\begin{proof}
	There exist eventually continuous maps $f:(X,Y) \to (CY,Y)$ and $g:(CY,Y) \to (X,Y)$ such that $f|_Y=g|_Y=1_Y$~\cite[Theorem 1.23]{cp}. By Lemma~\ref{equal}, $\mathscr{U}_f$ is an equivalence of categories with inverse $\mathscr{U}_g$. 
	
	Let $h:(X,Y) \to (X',Y')$ be an eventually continuous map. Define $f':(X',Y') \to (CY',Y')$ as in~\cite[Theorem 1.23]{cp}. Then $(Ch \circ f)|_Y = (f' \circ h)|_Y$. So by Lemma~\ref{equal}, $\mathscr{U}_{Ch} \circ \mathscr{U}_f$ and $\mathscr{U}_{f'} \circ \mathscr{U}_h$ are naturally equivalent.
\end{proof}

\begin{lem}\label{lem3}
	$C^*_r\mathcal{B}(X,Y;\mathbb{C})^W =C^*_r\mathcal{B}(X-(Y-W),W;\mathbb{C})^W$.
\end{lem}

\begin{proof}
	The two categories have precisely the same objects. Consider the functor from  $C_r^*\mathcal{B}(X,Y;\mathbb{C})^W$ to $C_r^* \mathcal{B}(X-(Y-W),W;\mathbb{C})^W$ induced by forgetting control along $Y-W$. To show that the functor is one-to-one on morphism sets, consider a morphism $T$ in $C_r^*\mathcal{B}(X,Y;\mathbb{C})$ whose image in $C_r^*\mathcal{B}(X-(Y-W),W;\mathbb{C})^W$ is zero. This means that $T$ can be approximated arbitrarily closely by continuously controlled bounded operators that factor through objects that are zero on a neighborhood of $W$. Therefore $T$ is identified with the zero morphism in $C_r^*\mathcal{B}(X,Y;\mathbb{C})^W$.
	
	To show that the functor is onto on morphism sets, let $T$ represent a morphism in $C_r^*\mathcal{B}(X-(Y-W),W;\mathbb{C})^W$. By~\cite[Proof of Lemma 1.32]{cp}, any morphism in $C_r^*\mathcal{B}(X-(Y-W),W;\mathbb{C})^W$ that can be represented by a continuously controlled operator is in the image of the functor. Since $T$ can be approximated arbitrarily well by continuously controlled operators, $T$ is related to a limit of morphisms in the image of the functor.
\end{proof}

\begin{lem}\label{lem3a}
	$C^*_r\mathcal{B}(CX,CY \cup X,p_X;\mathbb{C})^E = C^*_r\mathcal{B}(E \times (0,1],E \times 1;\mathbb{C})^{E \times 1}$.
\end{lem}

\begin{proof}
	By~\cite[Proof of Lemma 2.4]{cp}, the forgetful functor from $\mathcal{B} (CX,CY \cup X;\mathbb{C})^E$ to $\mathcal{B} (CX,CY \cup X,p_X;\mathbb{C})^E$ is an equivalence of categories. This equivalence extends to one between $C_r^*\mathcal{B} (CX,CY \cup X;\mathbb{C})^E$ and $C_r^*\mathcal{B} (CX,CY \cup X,p_X;\mathbb{C})^E$ as in the proof of Lemma~\ref{lem3}. Finally, $C_r^*\mathcal{B} (CX,CY \cup X;\mathbb{C})^E = C_r^*\mathcal{B} (E \times (0,1],E \times 1;\mathbb{C})^{E \times 1}$ by Lemma~\ref{lem3}.
\end{proof}

\begin{lem}\label{lem3b}
	Let $\Gamma$ be a group acting on $X$ such that $Y$ is $\Gamma$-invariant. Then
\[ \mathbb{K}^{Top}(C^*_r\mathcal{B}^{\Gamma}(CX,CY \cup X,p_X;\mathbb{C})) \simeq \mathbb{K}^{Top}(C^*_r\mathcal{B}^{\Gamma}(CX,CY \cup X,p_X;\mathbb{C})^E). \]
\end{lem}

\begin{proof}
	The Karoubi filtration $C^*_r\mathcal{B}^{\Gamma}(CX,CY \cup X,p_X;\mathbb{C})_{CY} \to C^*_r\mathcal{B}^{\Gamma}(CX,CY \cup X,p_X;\mathbb{C}) \to C^*_r\mathcal{B}^{\Gamma}(CX,CY \cup X,p_X;\mathbb{C})^E$ induces a fibration of spectra after applying $\mathbb{K}^{Top}$ and therefore produces a long exact sequence on homotopy groups. As with~\cite[Lemma 2.5]{cp}, the lemma is proven by an Eilenberg swindle. That is, we construct an endofunctor $\Sigma^{\infty}$ on $C^*_r\mathcal{B}^{\Gamma}(CX,CY \cup X,p_X;\mathbb{C})_{CY}$, such that $1 \oplus \Sigma^{\infty}$ is naturally equivalent to $\Sigma^{\infty}$. A category that admits such an endofunctor is called {\it flasque}. This implies that the spectrum $\mathbb{K}^{Top} (C^*_r\mathcal{B}^{\Gamma}(CX,CY \cup X,p_X;\mathbb{C})_{CY})$ is contractible.
	
	Choose a continuous function $s:X\to [1,\infty)$ such that $s(x)>1$ if $x\in E$ and $s(y)=1$ for every $y\in Y$. Let $H$ be an object in $C^*_r\mathcal{B}^{\Gamma}(CX,CY \cup X,p_X;\mathbb{C})_{CY}$. Let $S$ be the endofunctor on $\mathcal{B}^{\Gamma}(CX,CY \cup X,p_X;\mathbb{C})_{CY}$ defined by $S(H)_{(x,t)}=H_{(x,s(x)\cdot t)}$, where we set $H_{(x,s(x)\cdot t)}=0$ if $s(x)\cdot t\geq1$. Since $S$ is a unitary functor, it extends to an endofunctor $\mathscr{S}$ on $C^*_r\mathcal{B}^{\Gamma}(CX,CY \cup X,p_X;\mathbb{C})_{CY}$. Note that the unitary operator associated to $S$, $U_H:H \to S(H)$, is an isomorphism in $\mathcal{B}^{\Gamma}(CX,CY \cup X,p_X;\mathbb{C})_{CY}$.  Define $\Sigma^{\infty}(H)=\bigoplus_{n \geq 1}\mathscr{S}^n(H)$. If $T$ is a morphism in $C^*_r\mathcal{B}^{\Gamma}(CX,CY \cup X,p_X;\mathbb{C})_{CY}$, then $\Sigma^{\infty}(T)$ is the diagonal operator $\bigoplus_{n \geq 1}\mathscr{S}^n(T)$. Since the norms of the operators $\mathscr{S}^n(T)$ are uniformly bounded by $||T||$, $\Sigma^{\infty}(T)$ is bounded. Since objects are zero in a neighborhood of $E$, $\Sigma^{\infty}$ is a well-defined endofunctor. The natural equivalence, $\eta:1 \oplus \Sigma^{\infty}\to \Sigma^{\infty}$, is given by the diagonal unitary operator $\eta(H)=\bigoplus_{n \geq 0}U_{\mathscr{S}^n(H)}$.
\end{proof}

	It is worth noting that $p_X$-control plays an important role when $\Gamma$ is infinite. If we required continuous control everywhere, then equivariant morphisms could not have a non-zero component between $(p,s)$ and $(q,t)$ when $s\neq t$. The reason for this is that non-zero components of this type would be translated by the group to points close to $Y\times (0,1)$ contradicting control. The natural equivalence, however, requires non-zero components between such points.

\begin{coro}\label{isom}
	Let $\Gamma$ be a group acting on $X$ such that $Y$ is $\Gamma$-invariant. Then
\[ \mathbb{K}^{Top}(C^*_r\mathcal{B}^{\Gamma}(CX,CY \cup X,p_X;\mathbb{C})) \simeq \mathbb{K}^{Top}(C^*_r\mathcal{B}^{\Gamma}(E \times (0,1],E \times 1;\mathbb{C})^{E \times 1}). \]
\end{coro}

\begin{proof}
	Immediate from Lemmas~\ref{lem3a} and~\ref{lem3b}.
\end{proof}

	In~\cite{hp}, the {\it equivariant} continuously controlled category $\mathcal{B}_{\Gamma}(X,Y;R)$ is used to build equivariant homology theories. This category has a slightly different definition from the one used here. However, they agree when $\Gamma$ acts on $E=X-Y$ with finite isotropy. It was proved in~\cite[Theorem 7.5]{hp} that if $E$ is a cocompact $\Gamma$-space with finite isotropy, then
\[ \pi_i\big(\mathbb{K}^{Top}(C^*_r\mathcal{B}_{\Gamma}(E \times (0,1],E \times 1;\mathbb{C})^{E \times 1}) \big)\cong KK_{i-1}^{\Gamma}(C_0(E);\mathbb{C}). \]
Therefore, Corollary~\ref{isom} implies
\[ \pi_i\big(\mathbb{K}^{Top}(C^*_r\mathcal{B}^{\Gamma}(CX,CY \cup X,p_X;\mathbb{C}))\big)\cong KK_{i-1}^{\Gamma}(C_0(E);\mathbb{C}). \]
Given the following theorem, we see that the Baum-Connes assembly map is equivalent to 
\[ \Omega\mathbb{K}^{Top}(C^*_r\mathcal{B}^{\Gamma}(CX,CY \cup X,p_X;\mathbb{C})) \to \Omega\mathbb{K}^{Top}(C^*_r\mathcal{B}^{\Gamma}(\Sigma X,\Sigma Y,p_X;\mathbb{C})). \]

\begin{theo}\label{rightside}
	Assume that $\Gamma$, $\mathcal{E}$, $X$ and $Y$ satisfy the conditions of the main theorem. Then
\[ \Omega\mathbb{K}^{Top}(C^*_r\mathcal{B}^{\Gamma}(\Sigma X,\Sigma Y,p_X;\mathbb{C})) \simeq \mathbb{K}^{Top}(C^*_r\Gamma). \]
\end{theo}

\begin{proof}
	The assumption that the compact subsets of $\mathcal{E}$ shrink at infinity implies that the $p_X$-control condition is automatically satisfied~\cite[Proof of Lemma 2.3]{cp}. Thus, $C^*_r\mathcal{B}^{\Gamma}(\Sigma X,\Sigma Y,p_X;\mathbb{C})=C^*_r\mathcal{B}^{\Gamma}(\Sigma \mathcal{E},\{ 0,1 \};\mathbb{C})$. The quotient map $q:\Sigma \mathcal{E}\to \Sigma (\mathcal{E} / \Gamma)$, which is eventually continuous, induces an equivalence of categories
\[ \mathcal{B}^{\Gamma}(\Sigma \mathcal{E},\{ 0,1 \};\mathbb{C}) \to \mathcal{B}(\Sigma (\mathcal{E} / \Gamma),\{ 0,1 \};\mathbb{C}\Gamma). \]
To show this, we define the inverse functor $V$ as follows. Let $H$ be an object of $\mathcal{B}(\Sigma (\mathcal{E} / \Gamma),\{ 0,1 \};\mathbb{C}\Gamma)$. Recall that for every $z\in \Sigma (\mathcal{E} / \Gamma)$,  $H_z$ is a finitely generated free $\mathbb{C}\Gamma$-module. Therefore $H_z=\bigoplus_{g\in \Gamma}H_{(z,g)}$, where $H_{(z,g)}$ is a complex vector space whose dimension is equal to the $\mathbb{C}\Gamma$-dimension of $H_z$. Define $V(H)_{x}=\bigoplus_{g\in \Gamma_x}H_{(q(x),g)}$, where $\Gamma_x$ denotes the stabilizer subgroup of $x$. The fact that $V$ is a functor depends on the assumption that $\mathcal{E}$ is a finite $\Gamma$-CW complex. Furthermore, $V$ is a unitary functor since there is a corresponding unitary operator for each object. Thus the equivalence extends to its $C^*_r$-counterpart:
\[ C^*_r\mathcal{B}^{\Gamma}(\Sigma X,\Sigma Y,p_X;\mathbb{C}) \cong C^*_r\mathcal{B}(\Sigma (\mathcal{E} / \Gamma),\{ 0,1 \};\mathbb{C}\Gamma). \]
By Proposition~\ref{cone}, $C^*_r\mathcal{B}(\Sigma (\mathcal{E} / \Gamma),\{ 0,1 \};\mathbb{C}\Gamma)\cong C^*_r\mathcal{B}([0,1],\{ 0,1 \};\mathbb{C}\Gamma)$. To finish the proof, consider the commutative diagram
\[ \xymatrix{
	C^*_r\mathcal{B}_{\emptyset} \ar[d] \ar[r] & C^*_r\mathcal{B}_{ \{ 1 \} } \ar[d] \ar[r] & C^*_r\mathcal{B}_{ \{ 1 \} }^{ \{ 0,1 \} } \ar[d] \\
	C^*_r\mathcal{B}_{ \{ 0 \} } \ar[r] & C^*_r\mathcal{B} \ar[r] & C^*_r\mathcal{B}^{ \{ 1 \} } } \]
in which $\mathcal{B}=\mathcal{B}([0,1], \{ 0,1 \};\mathbb{C}\Gamma)$. Note that the rightmost vertical map is an equivalence of categories and that both $C^*_r\mathcal{B}_{ \{ 0 \} }$ and $C^*_r\mathcal{B}_{ \{ 1 \} }$ are flasque. Each row of the diagram is a Karoubi filtration inducing a fibration of spectra after applying $\mathbb{K}^{Top}$ and hence a long exact sequence of homotopy groups. Therefore, $\mathbb{K}^{Top}(C^*_r\mathcal{B}) \simeq \Sigma\mathbb{K}^{Top}(C^*_r\mathcal{B}_{\emptyset})$. The local finiteness of objects implies $C^*_r\mathcal{B}_{\emptyset} \cong C^*_r\mathcal{B}(\bullet, \emptyset;\mathbb{C}\Gamma)$. Thus, $\mathbb{K}^{Top}(C^*_r\mathcal{B}) \simeq \Sigma\mathbb{K}^{Top}(C^*_r\Gamma)$.
\end{proof}

\subsection{The Conner Conjecture and Steenrod Homology}\label{conner}

	Recall that a {\it reduced Steenrod homology theory}, $h$, is a functor from the category of compact metrizable spaces and continuous maps to the category of graded abelian groups satisfying
\begin{enumerate}
	\item[(i)] $h$ is homotopy invariant;
	\item[(ii)] $h(\bullet)=0$;
	\item[(iii)] given any closed subset $A$ of $X$, there is a natural transformation, $\partial :h_n(X / A) \to h_{n-1}(A)$, fitting into a long exact sequence
\[ \dotsb \to h_n(A) \to h_n(X) \to h_n(X / A) \to h_{n-1}(A) \to \dotsb ; \]
	\item[(iv)] if $\bigvee^{}_{} X_i$ denotes a compact metric space that is a countable union of metric spaces along a single point, then the projection maps $p_i : \bigvee^{}_{} X_i \to X_i$ induce an isomorphism $h_*(\bigvee^{}_{} X_i) \to \prod^{}_{} h_*(X_i)$.
\end{enumerate}
Given any generalized homology theory, there is a unique Steenrod homology extension. Existence of such extensions was proved by Kahn, Kaminker and Schochet, and Edwards and Hastings~\cite{kahn,edwards}. Uniqueness was proved by Milnor \cite{milnor}.

	A functor $k$ from the category of compact metrizable spaces to the category of spectra is called a {\it reduced Steenrod functor} if it satisfies the following conditions.
\begin{enumerate}
	\item[(i)] The spectrum $k(CX)$ is contractible.
	\item[(ii)] If $A \subset X$ is closed, then $k(A) \to k(X) \to k(X / A)$ is a fibration (up to natural weak homotopy equivalence).
	\item[(iii)] The projection maps $p_i : \bigvee^{}_{} X_i \to X_i$ induce a weak homotopy equivalence $k(\bigvee^{}_{} X_i) \to \prod^{}_{} k(X_i)$.
\end{enumerate}

\begin{theo}\label{sthom}\cite[Proposition 12.1]{higroeped}
        Let $\Gamma$ be a group. The functor
\[ \mathbb{K}^{Top}(C^*_r\mathcal{B}(C(-),-;\mathbb{C}\Gamma)) \]
is a reduced Steenrod functor whose value on $S^0$ is $\Sigma \mathbb{K}^{Top}(C^*_r\Gamma)$.
\end{theo}

	Given the work done thus far, it is easy to see that the proof of~\cite[Theorem 2.13]{cp} extends to give the following result. This proposition plays a key role in the proof of Theorem~\ref{filtration}.

\begin{prop}\label{key}
	Let $\Gamma$ be a group, $X$ a compact metrizable space, and $Y$ a closed nowhere dense subset. If the reduced Steenrod homology of $X$ (from Theorem~\ref{sthom}) is trivial, then
\[ \mathbb{K}^{Top}(C^*_r\mathcal{B}(CX,CY \cup X,p_X;\mathbb{C}\Gamma)) \simeq \mathbb{K}^{Top}(C^*_r\mathcal{B}(\Sigma X,\Sigma Y,p_X;\mathbb{C}\Gamma)). \]
\end{prop}

	In the proof of Theorem~\ref{filtration} we will be considering spaces with a finite group action. In order to use Proposition~\ref{key} on the corresponding quotient spaces, we need to know that the reduced Steenrod homology (from Theorem~\ref{sthom}) of the quotient of a contractible compact metrizable space by a finite group is trivial. The first step toward proving this is~\cite[Theorem III.7.12]{bredon}, which states that if $G$ is a finite group acting on a compact Hausdorff space $X$ that has trivial reduced \v{C}ech cohomology, then the reduced \v{C}ech cohomology of $X / G$ is also trivial. This is a special case of the Conner Conjecture which states that if a compact Lie group $G$ acts on a space $X$, where $X$ is either paracompact of finite cohomological dimension with finitely many orbit types or compact Hausdorff, then $X / G$ has trivial reduced \v{C}ech cohomology if $X$ does \cite{conner}. The Conner Conjecture was proved by Oliver~\cite{oliver}. The final step is the following theorem. 

\begin{theo}\cite[Theorem 5.2]{me2}\label{geogh}
	If the reduced \v{C}ech cohomology of a compact metrizable space $X$ is trivial, then  every reduced Steenrod homology of $X$ is also trivial.
\end{theo}

	As mentioned above, the following corollary will be used in the proof of Theorem~\ref{filtration}.

\begin{coro}\label{zero}
	Let $G$ be a finite group acting on a compact metrizable space $X$. If the reduced \v{C}ech cohomology of $X$ is trivial, then the reduced Steenrod homology of $X / G$ (from Theorem~\ref{sthom}) is trivial.
\end{coro}

\section{The Main Theorem}

\begin{theorem}
	Let $\Gamma$ be a discrete group, and let $\mathcal{E}=E\Gamma(\mathfrak{f})$ be the universal space for $\Gamma$-actions with finite isotropy, where $\mathfrak{f}$ denotes the family of finite subgroups of $\Gamma$. Assume that $\mathcal{E}$ is a finite $\Gamma$-CW complex admitting a compactification, $X$, (i.e., $X$ is compact and $\mathcal{E}$ is an open dense subset) such that
	\begin{enumerate}
		\item[1.] the $\Gamma$-action extends to $X$;
		\item[2.] $X$ is metrizable;
		\item[3.] $X^G$ is contractible for every $G\in \mathfrak{f}$;
		\item[4.] ${\mathcal{E}}^G$ is dense in $X^G$ for every $G\in \mathfrak{f}$;
		\item[5.] compact subsets of $\mathcal{E}$ become small near $Y=X-\mathcal{E}$.  That is, for every compact subset $K \subset \mathcal{E}$ and for every neighborhood $U \subset X$ of $y \in Y$, there exists a neighborhood $V \subset X$ of $y$ such that $g \in \Gamma$ and $gK \cap V \neq \emptyset$ implies $gK \subset U$.
	\end{enumerate}
Then the Baum-Connes map, $KK_i^{\Gamma} (C_0(\mathcal{E});\mathbb{C}) \to K_i(C^*_r\Gamma)$, is a split injection.
\end{theorem}

	There are many important classes of groups satisfying the conditions of the main theorem. {\it Crystallographic groups}, which are discrete groups that act cocompactly on Euclidean $n$-space by isometries, satisfy these conditions. The desired compactification is obtained by adding an $(n-1)$-sphere at infinity. More generally, {\it virtually polycyclic groups} satisfy these conditions, since we can also take $\mathbb{R}^n$ for some $n$ to be our universal space~\cite{wilking}. Gromov's {\it word hyperbolic groups} also satisfy the conditions by taking a certain compactification of the Rips complex~\cite{bestvina, meintrup}. It is already known that the the Baum-Connes assembly map is a split injection for word hyperbolic groups. Higson proved that the Baum-Connes map is a split injection for countable groups that admit an amenable action on a compact Hausdorff space~\cite{higson}, and a word hyperbolic group acts amenably on its Gromov boundary~\cite{adams}.

\subsection{The Proof}

	As was shown in Section~3.1, the map
\[ \Omega\mathbb{K}^{Top}(C^*_r\mathcal{B}^{\Gamma}(CX,CY \cup X,p_X;\mathbb{C})) \to \Omega\mathbb{K}^{Top}(C^*_r\mathcal{B}^{\Gamma}(\Sigma X,\Sigma Y,p_X;\mathbb{C})) \]
produces the Baum-Connes map after taking homotopy groups. Although it is not a map of fixed spectra, it fits in to the following commutative diagram in which $\mathcal{B}(CX)=\mathcal{B}(CX,CY \cup X,p_X;\mathbb{C})$ and $\mathcal{B}(\Sigma X)=\mathcal{B}^{\Gamma}(\Sigma X,\Sigma Y,p_X;\mathbb{C})$.
\[ \xymatrix{
	\Omega\mathbb{K}^{Top}(C^*_r\mathcal{B}^{\Gamma}(CX)) \ar[d] \ar[r] & \Omega\mathbb{K}^{Top}(C^*_r\mathcal{B}^{\Gamma}(\Sigma X)) \ar[d] \\
	\Omega\mathbb{K}^{Top}(C^*_r\mathcal{B}(CX))^{\Gamma} \ar[d] \ar[r] & \Omega\mathbb{K}^{Top}(C^*_r\mathcal{B}(\Sigma X))^{\Gamma} \ar[d] \\
	\Omega\mathbb{K}^{Top}(C^*_r\mathcal{B}(CX)^{h \mathfrak{f} \Gamma} \ar[r] & \Omega\mathbb{K}^{Top}(C^*_r\mathcal{B}(\Sigma X))^{h \mathfrak{f} \Gamma} } \]
Therefore we are able to make use of homotopy fixed sets as discussed in Section~2.4. The main theorem is proven by showing that the composition
\[ \Omega\mathbb{K}^{Top}(C^*_r\mathcal{B}^{\Gamma}(CX))\to \Omega\mathbb{K}^{Top}(C^*_r\mathcal{B}(CX))^{\Gamma}\to \Omega\mathbb{K}^{Top}(C^*_r\mathcal{B}(CX))^{h\mathfrak{f}\Gamma}, \]
and the map
\[ \Omega\mathbb{K}^{Top}(C^*_r\mathcal{B}(CX))^{h\mathfrak{f}\Gamma}\to \Omega\mathbb{K}^{Top}(C^*_r\mathcal{B}(\Sigma X))^{h\mathfrak{f}\Gamma}, \]
are weak homotopy equivalences. In this section the first of these two equivalences is handled. The second is proven in the subsequent section.

	For the remainder of this paper, assume that $\Gamma$, $\mathcal{E}$, $X$, and $Y$ satisfy the conditions of the main theorem. For notational convenience, let $\mathcal{B}(E \times (0,1);\mathbb{C})=\mathcal{B}(E \times (0,1],E \times 1;\mathbb{C})$.

\begin{theo}\label{theo3}
	The spectrum $\mathbb{K}^{Top}(C^*_r\mathcal{B}^{\Gamma}(\mathcal{E} \times (0,1);\mathbb{C})^{\mathcal{E} \times 1})$ is weakly homotopy equivalent to $\mathbb{K}^{Top}(C^*_r\mathcal{B}(\mathcal{E} \times (0,1);\mathbb{C})^{\mathcal{E} \times 1})^{h \mathfrak{f} \Gamma}$.
\end{theo}

\begin{proof}
Note that $\mathcal{E}$ is assumed to be a finite $\Gamma$-CW complex, and proceed by induction on the $\Gamma$-cells in $\mathcal{E}$. Begin with the discrete space ${\Gamma}/H$ for some $H \in \mathfrak{f}$. The control condition implies that the components of a morphism near ${\Gamma}/H \times 1$ must be zero between points with different ${\Gamma}/H$ entries. Since we are taking germs at ${\Gamma}/H \times 1$, the category $\mathcal{B}({\Gamma}/H \times (0,1);\mathbb{C})^{{\Gamma}/H \times 1}$ is equivalent to the product category $\prod_{\Gamma /H}^{} \mathcal{B}((0,1);\mathbb{C})^{1}$. Thus
\[ C^*_r\mathcal{B}({\Gamma}/H \times (0,1);\mathbb{C})^{{\Gamma}/H \times 1} \cong \bigoplus_{\Gamma /H}^{} C^*_r\mathcal{B}((0,1);\mathbb{C})^{1}. \]
The $\Gamma$-action on these categories is the one described in Lemma~\ref{lemmab}. Therefore, using Lemma~\ref{finite},
\begin{align}
	&\big(C^*_r\mathcal{B}({\Gamma}/H \times (0,1);\mathbb{C})^{{\Gamma}/H \times 1}\big)^{\Gamma} & &\cong & &\Big(\bigoplus_{\Gamma /H}^{} C^*_r\mathcal{B}((0,1);\mathbb{C})^{1}\Big)^{\Gamma} \notag\\ 
	& & &\cong & &\big(C^*_r\mathcal{B}((0,1);\mathbb{C})^{1}\big)^H \notag\\
	& & &\cong & &C^*_r\mathcal{B}^H((0,1);\mathbb{C})^{1} \notag\\
	& & &\cong & &C^*_r\Big(\big(\bigoplus_{\Gamma /H}^{} \mathcal{B}((0,1);\mathbb{C})^{1}\big)^{\Gamma}\Big) \notag\\
	& & &\cong & &C^*_r\mathcal{B}^{\Gamma}({\Gamma}/H \times (0,1);\mathbb{C})^{{\Gamma}/H \times 1}. \notag
\end{align}
Taking fixed sets commutes with applying $\mathbb{K}^{Top}$, therefore
\[ \mathbb{K}^{Top}(C^*_r\mathcal{B}^{\Gamma}({\Gamma}/H \times (0,1);\mathbb{C})^{{\Gamma}/H \times 1}) \cong \mathbb{K}^{Top}(C^*_r\mathcal{B}({\Gamma}/H \times (0,1);\mathbb{C})^{{\Gamma}/H \times 1})^{\Gamma}. \]
We now need to show that $\mathbb{K}^{Top}(C^*_r\mathcal{B}({\Gamma}/H \times (0,1);\mathbb{C})^{{\Gamma}/H \times 1})^{\Gamma}$ is weakly homotopy equivalent to $\mathbb{K}^{Top}(C^*_r\mathcal{B}({\Gamma}/H \times (0,1);\mathbb{C})^{{\Gamma}/H \times 1})^{h \mathfrak{f} \Gamma}$.

	The projection maps induce a map
\[ \mathbb{K}^{Top}\Big(\bigoplus_{\Gamma /H}^{} C^*_r\mathcal{B}((0,1);\mathbb{C})^{1}\Big) \to \prod_{\Gamma /H}^{} \mathbb{K}^{Top}(C^*_r\mathcal{B}((0,1);\mathbb{C})^{1}), \]
that is $\Gamma$-equivariant and a weak homotopy equivalence~(see Section 2.2). Consider the following commutative diagram:
\[ \xymatrix{
	\mathbb{K}^{Top}\Big(\bigoplus_{\Gamma /H}^{} C^*_r\mathcal{B}((0,1);\mathbb{C})^{1}\Big)^{\Gamma} \ar[d]_a \ar[r]^b &  \Big(\prod_{\Gamma /H}^{} \mathbb{K}^{Top}(C^*_r\mathcal{B}((0,1);\mathbb{C})^{1})\Big)^{\Gamma} \ar[d]_c \\
	\mathbb{K}^{Top}\Big(\bigoplus_{\Gamma /H}^{} C^*_r\mathcal{B}((0,1);\mathbb{C})^{1}\Big)^{h\mathfrak{f} \Gamma} \ar[r]^d & \Big(\prod_{\Gamma /H}^{} \mathbb{K}^{Top}(C^*_r\mathcal{B}((0,1);\mathbb{C})^{1})\Big)^{h\mathfrak{f} \Gamma}. } \]
We want to show that $a$ is a weak homotopy equivalence.

		Let $G \leq \Gamma$ be given, and choose representatives $\gamma_j$, so that $G \backslash \Gamma /H=\{ G\gamma_jH \}$. Then 
\begin{align}
	&\mathbb{K}^{Top}\Big(\bigoplus_{\Gamma /H}^{} C^*_r\mathcal{B}((0,1);\mathbb{C})^{1}\Big)^G & &\cong & &\mathbb{K}^{Top}\Big(\Big(\bigoplus_{\Gamma /H}^{} C^*_r\mathcal{B}((0,1);\mathbb{C})^{1}\Big)^G\Big) \notag\\ 
	& & &\cong & &\mathbb{K}^{Top}\Big(\bigoplus_{G \backslash \Gamma /H} C^*_r\mathcal{B}^{\gamma_j^{-1}G\gamma_j \cap H}((0,1],1;\mathbb{C})^{1} \Big) \notag\\
	& & &\simeq & &\prod_{G \backslash \Gamma /H}\mathbb{K}^{Top}\Big(C^*_r\mathcal{B}^{\gamma_j^{-1}G\gamma_j \cap H}((0,1],1;\mathbb{C})^{1} \Big) \notag\\
	& & &\simeq & &\prod_{G \backslash \Gamma /H}\mathbb{K}^{Top}(C^*_r\mathcal{B}((0,1],1;\mathbb{C})^{1})^{\gamma_j^{-1}G\gamma_j \cap H} \notag\\
	& & &\cong & &\Big(\prod_{G \backslash \Gamma /H}\mathbb{K}^{Top}(C^*_r\mathcal{B}((0,1],1;\mathbb{C})\Big)^G \notag
\end{align}
by Lemma~\ref{finite} and the fact that $\mathbb{K}^{Top}$ commutes with taking fixed sets. The case $G=\Gamma$ proves that $b$ is a weak homotopy equivalence. By Lemma $\ref{lemmaa}$, $d$ is a weak homotopy equivalence, and finally, $c$ is a weak homotopy equivalence by Lemma $\ref{lemmab}$. This completes the base case of the induction.

	Now assume that the theorem holds for $N$ and that $E$ is obtained from $N$ by attaching a $\Gamma$-$n$-cell, $\Gamma / K \times \mathbb{D}^n$, for some $K \in \mathfrak{f}$. Since $\mathcal{B}^{\Gamma}(E \times (0,1);\mathbb{C})_{N \times 1}^{E \times 1}$ is equivalent to $\mathcal{B}^{\Gamma}(N \times (0,1);\mathbb{C})^{N \times 1}$, given by the projection functor, which is unitary,
\[ C^*_r\mathcal{B}^{\Gamma}(N \times (0,1);\mathbb{C})^{N \times 1} \to C^*_r\mathcal{B}^{\Gamma}(E \times (0,1);\mathbb{C})^{E \times 1} \to C^*_r\mathcal{B}^{\Gamma}(E \times (0,1);\mathbb{C})^{(E-N) \times 1} \]
is a Karoubi filtration and therefore induces a fibration of spectra after applying $\mathbb{K}^{Top}$. The same statement is true if $\Gamma$ is replaced by a finite subgroup $G$. For notational convenience, let $\mathscr{A}=\mathcal{B}(N \times (0,1);\mathbb{C})^{N \times 1}$, $\mathscr{B}=\mathcal{B}(E \times (0,1);\mathbb{C})^{E \times 1}$, and $\mathscr{C}=\mathcal{B}(E \times (0,1);\mathbb{C})^{(E-N) \times 1}$. 
	
	Consider the following commutative diagram:
\[ \xymatrix{
	\mathbb{K}^{Top}(C^*_r\mathscr{A}^{\Gamma}) \ar[d]_a \ar[r] & \mathbb{K}^{Top}(C^*_r\mathscr{B}^{\Gamma}) \ar[d]_b \ar[r] & \mathbb{K}^{Top}(C^*_r\mathscr{C}^{\Gamma}) \ar[d]_c \\
	\mathbb{K}^{Top}(C^*_r\mathscr{A})^{h \mathfrak{f} \Gamma} \ar[r] & \mathbb{K}^{Top}(C^*_r\mathscr{B})^{h \mathfrak{f} \Gamma} \ar[r] & \mathbb{K}^{Top}(C^*_r\mathscr{C})^{h \mathfrak{f} \Gamma}. } \]
By the induction hypothesis, $a$ is a weak homotopy equivalence. To complete the proof, we need to show that $b$ is a weak homotopy equivalence. The top row in the diagram is a fibration of spectra. By the Five Lemma, it suffices to prove that $c$ is a weak homotopy equivalence once we know that the bottom row is also a fibration of spectra. As mentioned above, $C^*_r\mathscr{A}^G \to C^*_r\mathscr{B}^G \to C^*_r\mathscr{C}^G$ induces a fibration after applying $\mathbb{K}^{Top}$, for every $G \in \mathfrak{f}$. Thus, by Lemma~\ref{finite},
\[ \mathbb{K}^{Top}(C^*_r\mathscr{A})^G \to \mathbb{K}^{Top}(C^*_r\mathscr{B})^G \to \mathbb{K}^{Top}(C^*_r\mathscr{C})^G \]
is a fibration for every $G\in \mathfrak{f}$. Let $F$ denote the homotopy fiber of $\mathbb{K}^{Top}(C^*_r\mathscr{B}) \to \mathbb{K}^{Top}(C^*_r\mathscr{C})$. Taking homotopy fixed sets and homotopy fibers are commutative operations since both are inverse limits. Therefore $F^{h \mathfrak{f} \Gamma} \to \mathbb{K}^{Top}(C^*_r\mathscr{B})^{h \mathfrak{f} \Gamma} \to \mathbb{K}^{Top}(C^*_r\mathscr{C})^{h \mathfrak{f} \Gamma}$ is a fibration of spectra. Taking fixed sets is also an inverse limit, therefore it too commutes with homotopy fibers. Hence, $F^G \to \mathbb{K}^{Top}(C^*_r\mathscr{B})^G \to \mathbb{K}^{Top}(C^*_r\mathscr{C})^G$ is a fibration for every $G \in \mathfrak{f}$. Therefore $F^G \simeq \mathbb{K}^{Top}(C^*_r\mathscr{A})^G$ for every $G \in \mathfrak{f}$. By Lemma $\ref{lemmaa}$, $F^{h \mathfrak{f} \Gamma} \simeq \mathbb{K}^{Top}(C^*_r\mathscr{A})^{h \mathfrak{f} \Gamma}$. Therefore the bottom row in the above diagram is a fibration.
	
	Since $E-N= \Gamma / K \times \mathring{e^n}$, where $\mathring{e^n}$ is an open $n$-cell, the category $\mathcal{B} (E \times (0,1);\mathbb{C})^{(E-N) \times 1}$ is equivalent to the product category $\prod_{\Gamma /K}^{} \mathcal{B} (\mathring{e^n} \times (0,1);\mathbb{C})^{\mathring{e^n} \times 1}$. But this is entirely similar to the start of the induction. Therefore $c$ is a weak homotopy equivalence.
\end{proof}

\begin{coro}
	The spectrum $\mathbb{K}^{Top}(C^*_r\mathcal{B}^{\Gamma}(CX,CY \cup X,p_X;\mathbb{C}))$ is weakly homotopy equivalent to $\mathbb{K}^{Top}(C^*_r\mathcal{B}(CX,CY \cup X,p_X;\mathbb{C}))^{h \mathfrak{f} \Gamma}$.
\end{coro}

\begin{proof}
\begin{align}
	&\mathbb{K}^{Top}(C^*_r\mathcal{B}^{\Gamma}(CX,CY \cup X,p_X;\mathbb{C})) & &\simeq & &\mathbb{K}^{Top}(C^*_r\mathcal{B}^{\Gamma}(CX,CY \cup X,p_X;\mathbb{C})^\mathcal{E}) \notag\\ 
	& & &\cong & &\mathbb{K}^{Top}(C^*_r\mathcal{B}^{\Gamma}(\mathcal{E} \times (0,1);\mathbb{C})^{\mathcal{E} \times 1}) \notag\\
	& & &\simeq & &\mathbb{K}^{Top}(C^*_r\mathcal{B}(\mathcal{E} \times (0,1);\mathbb{C})^{\mathcal{E} \times 1})^{h \mathfrak{f} \Gamma} \notag\\
	& & &\cong & &\mathbb{K}^{Top}(C^*_r\mathcal{B}(CX,CY \cup X,p_X;\mathbb{C})^\mathcal{E})^{h \mathfrak{f} \Gamma} \notag\\
	& & &\simeq & &\mathbb{K}^{Top}(C^*_r\mathcal{B}(CX,CY \cup X,p_X;\mathbb{C}))^{h \mathfrak{f} \Gamma}. \notag
\end{align}
	
	The first three equivalences follow from Lemma~\ref{lem3b}, Lemma~\ref{lem3a}, and Theorem~\ref{theo3}, respectively. Lemma~\ref{lem3b} also holds when $\Gamma$ is replaced by any of its finite subgroups. Therefore, by Lemma~\ref{finite}, $\mathbb{K}^{Top}(C^*_r\mathcal{B}(CX,CY \cup X,p_X;\mathbb{C}))^G \cong \mathbb{K}^{Top}(C^*_r\mathcal{B}(CX,CY \cup X,p_X;\mathbb{C})^\mathcal{E})^G$ for every $G\in \mathfrak{f}$. Now Lemma~\ref{lemmaa} implies that $\mathbb{K}^{Top}(C^*_r\mathcal{B}(CX,CY \cup X,p_X;\mathbb{C})^{h \mathfrak{f} \Gamma} \simeq \mathbb{K}^{Top}(C^*_r\mathcal{B}(CX,CY \cup X,p_X;\mathbb{C})^\mathcal{E})^{h \mathfrak{f} \Gamma}$. Finally, using a similar argument with Lemma~\ref{lem3a}, $\mathbb{K}^{Top}(C^*_r\mathcal{B}(\mathcal{E} \times (0,1);\mathbb{C})^{\mathcal{E} \times 1})^{h \mathfrak{f} \Gamma} \cong \mathbb{K}^{Top}(C^*_r\mathcal{B}(CX,CY \cup X,p_X;\mathbb{C})^\mathcal{E})^{h \mathfrak{f} \Gamma}$.
\end{proof}

\subsection{Filtering by Conjugacy Classes of Fixed Sets}\label{filtering}

	Now we complete the proof of the main theorem by showing that
\[ \Omega\mathbb{K}^{Top}(C^*_r\mathcal{B}(CX,CY \cup X,p_X;\mathbb{C}))^{h\mathfrak{f}\Gamma} \simeq \Omega\mathbb{K}^{Top}(C^*_r\mathcal{B}(\Sigma X,\Sigma Y,p_X;\mathbb{C}))^{h\mathfrak{f}\Gamma}. \]
By Lemmas~\ref{finite} and~\ref{lemmaa}, it suffices to prove the following.
	
\begin{theo}\label{filtration}
	For every $G \in \mathfrak{f}$,
\[ \mathbb{K}^{Top}(C^*_r\mathcal{B}^G(CX,CY \cup X,p_X;\mathbb{C})) \simeq \mathbb{K}^{Top}(C^*_r\mathcal{B}^G(\Sigma X,\Sigma Y,p_X;\mathbb{C})). \]
\end{theo}

	Given $G \in \mathfrak{f}$, consider the subgroup lattice for $G$, and let $H\leq G$. Define the distance from $G$ to $H$, ${\rm dist}(H)$, to be the maximum number of steps needed to reach $H$ from $G$ on the subgroup lattice. Notice that ${\rm dist}(gHg^{-1})={\rm dist}(H)$ for each $g\in G$.

	Let $n={\rm dist}(1)$, and let $l_i$ be the number of conjugacy classes of subgroups with distance $i$ from $G$. For each $i$, $1\leq i\leq n-1$, choose a representative, $H_{i,j}$, $1\leq j\leq l_i$, for each of the conjugacy classes with distance $i$ from $G$. Order the $H_{i,j}$'s using the dictionary order on the indexing set. Now re-index the sequence of subgroups according to the ordering so that
\[H_1=H_{1,1} \ \ , \ \ H_2=H_{1,2} \ \ , \ ... \ , \ \ H_m=H_{n-1,l_{n-1}}. \]
For convenience define $H_0=G$ and $H_{m+1}=1$.

	For each $j$, $0\leq j\leq m+1$, define $C_j=\bigcup_{g \in G} X^{gH_jg^{-1}}$, which will often be referred to as a {\it conjugacy class of fixed sets}.

	For each $k$, $0\leq k\leq m+1$, define $Z_k=\bigcup_{0 \leq j \leq k} C_j$. Also define $Y_k=Z_k \cap Y$. Since it is possible that $Y_k= \emptyset$, we define $C(\emptyset)=\{ 0 \}$ and $\Sigma (\emptyset)=\{ 0,1 \}$ for the convenience of notation.
	
	Notice that $Z_k$ is contractible for every $k$, $0\leq k\leq m+1$. Also notice that since $gX^H=X^{gHg^{-1}}$, each of the conjugacy classes of fixed sets is $G$-invariant. Therefore $Z_k$ is $G$-invariant for every $k$, $0\leq k\leq m+1$.

	To simplify the notation set
\begin{align}
	&\mathcal{B}^G_{\mathbb{C}}(C(Z_k))= \mathcal{B}^G(C(Z_k),C(Y_k) \cup Z_k,p_{Z_k};\mathbb{C}); \notag\\ 
	&\mathcal{B}^G_{\mathbb{C}}(\Sigma (Z_k)) = \mathcal{B}^G(\Sigma (Z_k),\Sigma (Y_k),p_{Z_k};\mathbb{C}); \notag\\
	&\mathcal{B}^G_{\mathbb{C}}(C(Z_k);\mathbb{C})_{Z_{k-1}} = \mathcal{B}^G_{\mathbb{C}}(C(Z_k))_{C(Y_{k-1}) \cup Z_{k-1}}; \notag\\
        &\mathcal{B}^G_{\mathbb{C}}(\Sigma (Z_k))_{Z_{k-1}} = \mathcal{B}^G_{\mathbb{C}}(\Sigma (Z_k))_{\Sigma (Y_{k-1})}; \notag\\
	&\mathcal{B}^G_{\mathbb{C}}(C(Z_k))^{>Z_{k-1}} = \mathcal{B}^G_{\mathbb{C}}(C(Z_k))^{(C(Y_k) \cup Z_k)-(C(Y_{k-1}) \cup Z_{k-1})}; \notag\\
	&\mathcal{B}^G_{\mathbb{C}}(\Sigma (Z_k))^{>Z_{k-1}} = \mathcal{B}^G_{\mathbb{C}}(\Sigma (Z_k))^{\Sigma(Y_k)-\Sigma (Y_{k-1})} . \notag
\end{align}

	Theorem $\ref{filtration}$ is proved by induction on the chain \[X^G=Z_0 \subseteq Z_1 \subseteq \dddot{} \subseteq Z_m \subseteq Z_{m+1}=X. \]
As in~\cite{me2}, this is accomplished with the following three lemmas.

\begin{lem}\label{lemma1}
	Let $1\leq k\leq m+1$. Then there are involution preserving equivalences
	\begin{enumerate}
		\item[i.] $C^*_r\mathcal{B}^G_{\mathbb{C}}(C(Z_{k}))_{Z_{k-1}} \cong C^*_r\mathcal{B}^G_{\mathbb{C}}(C(Z_{k-1}))$, and
		\item[ii.] $C^*_r\mathcal{B}^G_{\mathbb{C}}(\Sigma (Z_{k}))_{Z_{k-1}} \cong C^*_r\mathcal{B}^G_{\mathbb{C}}(\Sigma (Z_{k-1}))$.
	\end{enumerate}
\end{lem}

\begin{proof}
	The proof of this lemma is given in~\cite[Lemma 3.3]{me3}, and is included here for the reader's convenience. Consider part (i). The equivalence is proven by constructing an equivariant function $f$, from $C(Z_k)$ into itself, that is not quite an eventually continuous map, but is good enough to induce a functor from $C^*_r\mathcal{B}^G_{\mathbb{C}}(C(Z_{k}))_{Z_{k-1}}$ to $C^*_r\mathcal{B}^G_{\mathbb{C}}(C(Z_{k-1}))$ that is an inverse to the inclusion functor up to natural equivalence. The only difficulty in defining such a function appears when there are sequences consisting of points outside of $C(Z_{k-1})$ that converge to points in $C(Y_{k-1}) \cup Z_{k-1}$. If this does not happen, then every point not in $C(Z_{k-1})$ is greater than a fixed distance away from $C(Y_{k-1}) \cup Z_{k-1}$. Since objects are zero on a neighborhood of $C(Y_k) \cup Z_k-(C(Y_{k-1}) \cup Z_{k-1})$, there will only be finitely many points in $C(Z_k)-C(Z_{k-1})$ for which an object is non-zero. Thus, we could define $f$ to be the identity on $C(Y_k)\cup Z_k \cup C(Y_{k-1})$ and to send every other point to a chosen point in $\mathcal{E}^G \times (0,1)$.
	
	Assume that such sequences exist. Choose a representative $x$, for each orbit not contained in $C(Z_{k-1})$, such that $G_x=H_k$. Choose $y_x\in (C(Y_{k-1}) \cup Z_{k-1})^{H_k}$ such that $d(x,y_x)=d(x,(C(Y_{k-1}) \cup Z_{k-1})^{H_k})$. This is possible by the assumption that there exist sequences consisting of points outside of $C(Z_{k-1})$ converging to points in $C(Y_{k-1}) \cup Z_{k-1}$. Choose $z_x\in (Z_{k-1}-Y_{k-1})^{H_k}\times (0,1)$ such that $d(z_x,y_x)\leq d(x, C(Y_{k-1}) \cup Z_{k-1})$. This choice is possible by condition 4 of the main theorem. Define $f:(C(Z_k),C(Y_k\cup Z_k)) \to (C(Z_k),C(Y_k\cup Z_k))$ by
\[ f(a)=\left\{
\begin{array}{cl}
	gz_x & {\rm if} \ a=gx \\
	a & {\rm if} \ a \in C(Y_k)\cup Z_k \cup C(Y_{k-1}).
\end{array}
\right.\]
This is well-defined since $g_1x=g_2x$ if and only if $g_2^{-1}g_1\in H_k$. Thus, $f$ is $G$-equivariant. As mentioned before, $f$ is not eventually continuous. More precisely, the second and third conditions of eventually continuous maps are not satisfied by $f$. The second condition ensures that the induced functor sends objects to objects. To replace this condition, we show that the inverse image under $f$ of a compact set in $C(Z_{k-1})-(C(Y_{k-1}) \cup Z_{k-1})$ contains only finitely many points over which an object can be non-zero. We do this as follows.

	Let $w\in C(Z_{k-1})-(C(Y_{k-1}) \cup Z_{k-1})$ and $d>0$ such that $0<d<d(w,C(Y_{k-1}) \cup Z_{k-1})$ be given. Let
\[ U=\big\{ z\in C(Z_k) \, \big| \, d(z,C(Y_{k-1}) \cup Z_{k-1})<d(w,C(Y_{k-1}) \cup Z_{k-1}-d) \big\}. \]
If $B_d(w)$ is the open ball in $C(Z_{k-1})$ about $w$ of radius $d$, then $f^{-1}(\overline{B_d(w)})\subseteq C(Z_k)-U$. To see why, let $a\in f^{-1}(\overline{B_d(w)})$. If $a\in C(Z_{k-1})$ then $a=f(a)\in \overline{B_d(w)}$. If $a\notin C(Z_{k-1})$ then $a=gx$, where $x$ is the chosen representative, and $d(w,C(Y_{k-1}) \cup Z_{k-1})-d \leq d(gz_x,C(Y_{k-1}) \cup Z_{k-1})=d(z_x, C(Y_{k-1}) \cup Z_{k-1}) \leq d(z_x, C(Y_{k-1}^{H_k}) \cup Z_{k-1}^{H_k}) \leq d(z_x,y_x) \leq d(x, C(Y_{k-1}) \cup Z_{k-1})=d(a, C(Y_{k-1}) \cup Z_{k-1})$. Given an object $H$, there is a neighborhood $V\subseteq C(Z_k)$ of $C(Y_k) \cup Z_k-(C(Y_{k-1}) \cup Z_{k-1})$ on which $H$ is zero. Therefore,
\[ \big\{a\in \overline{f^{-1}(\overline{B_d(w)})} \, \big| \, H_a \neq 0 \big\} \subset (C(Z_k)-U)\cap (C(Z_k)-V). \]
Since $(C(Z_k)-U)\cap (C(Z_k)-V) \subseteq C(Z_k)-(C(Y_k \cup Z_k)$ is compact, this set is finite.

	The third condition on eventually continuous maps is that they are continuous on the boundary. This condition guarantees that the image under the induced functor of a continuously controlled morphism is continuously controlled.  Since objects are zero on a neighborhood of $C(Y_k) \cup Z_k-(C(Y_{k-1}) \cup Z_{k-1})$, it suffices to show that $f$ is continuous on $C(Y_{k-1}) \cup Z_{k-1}$. Let $\{ x_n \}$ be a sequence in $C(Z_k)-C(Z_{k-1})$ converging to $y\in C(Y_{k-1}) \cup Z_{k-1}$. Since $G$ is a finite group, all but finitely many terms in the sequence must be contained in a subsequence of the form $\{ gx_m \}$ for a fixed $g\in G$, where each $x_k$ is the chosen representative of its orbit. Since $\{ x_m \}$ converges to $g^{-1}y$ and $G_{x_m}=H_k$, $g^{-1}y$ is fixed by $H_k$. Thus, $\{ y_{x_m} \}$ converges to $g^{-1}y$. Therefore, $\{ z_{x_m} \}$ converges to $g^{-1}y$ and $\{ f(gx_m)\}=\{ gz_{x_m} \}$ converges to $y=f(y)$. Since $f|_{C(Y_k) \cup Z_k}=1_{C(Y_k) \cup Z_k}$, $f$ induces a functor that is an inverse to the inclusion functor up to natural equivalence as in Lemma~\ref{equal}.
	
	This completes the proof of part (i). Part (ii) is proven the same way replacing cones with suspensions.
\end{proof}

\begin{lem}\label{lemma2}
	Let $1\leq k\leq m+1$. Then there are involution preserving equivalences
	\begin{enumerate}
		\item[i.] $C^*_r\mathcal{B}^G_{\mathbb{C}}(C(Z_k))^{>Z_{k-1}} \cong C^*_r\mathcal{B}_{\mathbb{C} H_k}(C(Z_k / G))^{>Z_{k-1} / G}$, and
		\item[ii.] $C^*_r\mathcal{B}^G_{\mathbb{C}}(\Sigma (Z_k))^{>Z_{k-1}} \cong C^*_r\mathcal{B}_{\mathbb{C} H_k}(\Sigma (Z_k / G))^{>Z_{k-1} / G}$.
	\end{enumerate}
\end{lem}

\begin{proof}
	For every $z\in C(Z_k / G)-C(Z_{k-1} / G)$ choose $x_z\in p^{-1}(z)$, where $p:C(Z_k) \to C(Z_k / G)$ denotes the quotient map, so that $G_{x_z}=H_k$. Suppose $H$ and $K$ are objects in $C^*_r\mathcal{B}^G_{\mathbb{C}}(C(Z_k))$ that are zero over points in $C(Z_{k-1})$, and that $T:H\to K$ is a continuously controlled morphism such that:
\[ T^{x'}_{y'}=\left\{
\begin{array}{cl}
	gT^x_yg^{-1} & {\rm if} \ x'=gx \ {\rm and} \ y'=gy \\
	0 & {\rm otherwise}
\end{array}
\right.\]
where $G_x=G_y=H_k$. If $T^{Gx}_{Gy}:H_{Gx}\to K_{Gy}$ is non-zero, where $H_{Gx}=\bigoplus_{x'\in Gx} H_{x'}$, then the non-zero components of $T^{Gx}_{Gy}$ are all unitarily equivalent to $T^x_y$. We shall call $T^x_y$ a {\it generating component}. Note that the equivariance of $T$ implies that $T^x_y:H_x\to K_y$ is $H_k$-equivariant. Define the object $\bar{H}$ in $C^*_r\mathcal{B}_{\mathbb{C} H_k}(C(Z_k / G))$ as follows: $\bar{H}_z=H_{x_z}$ if $z\notin C(Z_{k-1})$ and is zero otherwise. Define the morphism $\bar{T}:\bar{H}\to \bar{K}$ as follows: $\bar{T}^z_w=h^{-1}T^{gx_z}_{hx_w}g$, where $T^{gx_z}_{hx_w}$ is a generating component. This assignment is well-defined since $T$ is equivariant. Note that $\| \bar{T} \|=\| T \|$, so $\bar{T}$ will be bounded if and only if $T$ is bounded. This also implies that if $S=\lim_n T_n$ is a limit of such morphisms, then $\bar{S}=\lim_n \bar{T_n}$ is well-defined.

	The scenario described above tells the entire story. The reason is that since we are taking germs away from $C(Y_{k-1}) \cup Z_{k-1}$ (resp. $C(Y_{k-1} / G) \cup Z_{k-1} / G$), the components of continuously controlled morphisms need to become small. This implies that a morphism in $\mathcal{B}^G_{\mathbb{C}}(C(Z_{k}))^{>Z_{k-1}}$ (resp. $\mathcal{B}_{\mathbb{C}H_k}(C(Z_{k} / G))^{>Z_{k-1} / G}$) will have a representative that is zero on $C(Z_{k-1})$ (resp. $C(Z_{k-1} / G)$). It is therefore irrelevant what the objects over $C(Z_{k-1})$ and $C(Z_{k-1} / G)$ are. Since we are working with a finite group action, non-zero components of a continuously controlled morphism have the same isotropy, namely a conjugate of $H_k$. This also tells us that the equivariance of morphisms in $C^*_r\mathcal{B}^G_{\mathbb{C}}(C(Z_{k}))^{>Z_{k-1}}$ implies that there is only one choice when lifting a continuously controlled morphism from $C^*_r\mathcal{B}_{\mathbb{C}H_k}(C(Z_{k} / G))^{>Z_{k-1} / G}$. Furthermore, if two morphisms $S$ and $S'$, are limits of continuously controlled morphisms of the type described above and are identified in $C^*_r\mathcal{B}^G_{\mathbb{C}}(C(Z_{k}))^{>Z_{k-1}}$, then $\bar{S}$ will be related to $\bar{S'}$. Since every morphism can be identified with a morphism that is such a limit, the above work defines a functor from $C^*_r\mathcal{B}^G_{\mathbb{C}}(C(Z_k))^{>Z_{k-1}}$ to $C^*_r\mathcal{B}_{\mathbb{C}H_k}(C(Z_k / G))^{>Z_{k-1} / G}$ that is an equivalence of categories. The same argument proves the second part of this lemma, replacing cones with suspensions.
\end{proof}

\begin{lem}\label{lemma3}
	For each $k$, $1\leq k\leq m+1$,
\[ \mathbb{K}^{Top}\big(C^*_r\mathcal{B}_{\mathbb{C} H_k}(C(Z_k / G))^{>Z_{k-1} / G}\big) \simeq \mathbb{K}^{Top}\big(C^*_r\mathcal{B}_{\mathbb{C} H_k}(\Sigma (Z_k / G))^{>Z_{k-1} / G}\big). \]
\end{lem}

\begin{proof}
	The proof proceeds just as in~\cite[Lemma 7.5]{me2}. Consider the following commutative diagram:
\[ \xymatrix{
	C^*_r\mathcal{B}_{\mathbb{C}H_k}(C({Z_k} / G))_{Z_{k-1} / G} \ar[d]_a \ar[r] & C^*_r\mathcal{B}_{\mathbb{C}H_k}(C({Z_k} / G)) \ar[d]_b \ar[r] & C^*_r\mathcal{B}_{\mathbb{C}H_k}(C({Z_k} / G))^{>Z_{k-1} / G} \ar[d]_c \\
	C^*_r\mathcal{B}_{\mathbb{C}H_k}(\Sigma({Z_k} / {G}))_{Z_{k-1} / G} \ar[r] & C^*_r\mathcal{B}_{\mathbb{C}H_k}(\Sigma({Z_k} / {G})) \ar[r] & C^*_r\mathcal{B}_{\mathbb{C}H_k}(\Sigma({Z_k} / {G}))^{>Z_{k-1} / G}. } \]
Since $Z_k$ is contractible, $b$ induces a weak homotopy equivalence by Corollary~\ref{zero} and Proposition~\ref{key}. An unequivariant version of the proof of Lemma~\ref{lemma1} proves
\[ C^*_r\mathcal{B}_{\mathbb{C}H_k}(C({Z_k} / G))_{Z_{k-1} / G}\cong C^*_r\mathcal{B}_{\mathbb{C}H_k}(C(Z_{k-1} / G)) \]
and
\[  C^*_r\mathcal{B}_{\mathbb{C}H_k}(\Sigma({Z_k} / G))_{Z_{k-1} / G}\cong C^*_r\mathcal{B}_{\mathbb{C}H_k}(\Sigma(Z_{k-1} / G)).\]
Then Corollary~\ref{zero} and Proposition~\ref{key} apply again to show that $a$ induces a weak homotopy equivalence. Since each row in the diagram is a Karoubi filtration, the Five Lemma shows that $c$ also induces a weak homotopy equivalence.
\end{proof}

We are now ready to prove Theorem~\ref{filtration}.

\begin{proof}[Proof of Theorem~\ref{filtration}]
	Since $G$ acts trivially on $X^G$,
\[ C^*_r\mathcal{B}^G_{\mathbb{C}}(C(X^G)) \cong C^*_r\mathcal{B}_{\mathbb{C}G}(C(X^G)) \]
and
\[ C^*_r\mathcal{B}^G_{\mathbb{C}}(\Sigma(X^G)) \cong C^*_r\mathcal{B}_{\mathbb{C}G}(\Sigma(X^G)). \]
Since $X^G$ is contractible, Proposition~\ref{key} implies
\[ \mathbb{K}^{Top}(C^*_r\mathcal{B}_{\mathbb{C}G}(C(X^G))) \simeq \mathbb{K}^{Top}(C^*_r\mathcal{B}_{\mathbb{C}G}(\Sigma(X^G))). \]
Therefore, $\mathbb{K}^{Top}\big(C^*_r\mathcal{B}^G_{\mathbb{C}}(C(X^G))\big) \simeq \mathbb{K}^{Top}\big(C^*_r\mathcal{B}^G_{\mathbb{C}}(\Sigma(X^G))\big)$. This completes the base case of the induction. 
	
	Assume now that $\mathbb{K}^{Top}(C^*_r\mathcal{B}^G_{\mathbb{C}}(C(Z_{k-1}))) \simeq \mathbb{K}^{Top}(C^*_r\mathcal{B}^G_{\mathbb{C}}(\Sigma(Z_{k-1})))$. We want to show that $\mathbb{K}^{Top}(C^*_r\mathcal{B}^G_{\mathbb{C}}(C(Z_k))) \simeq \mathbb{K}^{Top}(C^*_r\mathcal{B}^G_{\mathbb{C}}(\Sigma(Z_k)))$.
Consider the following commutative diagram:
\[ \xymatrix{
	C^*_r\mathcal{B}^G_{\mathbb{C}}(C(Z_k))_{Z_{k-1}} \ar[d]_a \ar[r] & C^*_r\mathcal{B}^G_{\mathbb{C}}(C(Z_k)) \ar[d]_b \ar[r] & C^*_r\mathcal{B}^G_{\mathbb{C}}(C(Z_k))^{>Z_{k-1}} \ar[d]_c \\
	C^*_r\mathcal{B}^G_{\mathbb{C}}(\Sigma(Z_k))_{Z_{k-1}} \ar[r] & C^*_r\mathcal{B}^G_{\mathbb{C}}(\Sigma(Z_k)) \ar[r] & C^*_r\mathcal{B}^G_{\mathbb{C}}(\Sigma(Z_k))^{>Z_{k-1}}. } \]

	By Lemma $\ref{lemma1}$, Proposition $\ref{key}$, and the induction hypothesis, $a$ induces a weak homotopy equivalence. By Lemmas $\ref{lemma2}$ and $\ref{lemma3}$, $c$ induces a weak homotopy equivalence. Since each row in the diagram is a Karoubi filtration, the Five Lemma shows that $b$ also induces a weak homotopy equivalence.
 \end{proof}

	This completes the proof of the main theorem.

\end{document}